# W-SURFACES HAVING SOME PROPERTIES

Stylianos Stamatakis

**Abstract.** We investigate some characteristic properties of specific Weingarten surfaces in the Euclidean space $E^3$ using the nets of the lines of curvature resp. the asymptotic lines on both central surfaces of them.

**2000 Mathematics Subject Classification:** 53A05

**Keywords:** Surfaces in Euclidean space, Weingarten surfaces

**1.** A surface of *Weingarten*, or simply W-surface, is a surface in the three-dimensional space $E^3$ for which there is a functional relation between the principle curvatures $k_1$, $k_2$. A well-known characterization of the *W*-surfaces is:

*A surface $\mathfrak{F}$ is a W-surface if and only if the second fundamental forms $\mathrm{II}^{(i)}$, $i = 1, 2$, of the central surfaces $\mathfrak{F}^{(i)}$ of $\mathfrak{F}$ are proportional.*

In case that the central surfaces $\mathfrak{F}^{(i)}$ are of negative Gaussian curvature the later property means:

*The asymptotic lines on both central surfaces correspond to each other.*

This theorem is due to *A. Ribaucour* as well as the following one [1, p. 184]:

*The lines of curvature on both central surfaces correspond to each other if and only if $\mathfrak{F}$ is a W-surface satisfying $\dfrac{1}{k_1} - \dfrac{1}{k_2} = const.$*

We also note the property [1, p.180]:

*The nets on $\mathfrak{F}^{(i)}$, $i = 1, 2$, corresponding to the lines of curvature of $\mathfrak{F}$, are both conjugate.*

In this paper we consider the nets of the lines of curvature resp. the asymptotic lines on both central surfaces $\mathfrak{F}^{(i)}$, $i = 1, 2$, of a given surface $\mathfrak{F}$ and the corresponding nets on $\mathfrak{F}$ and we ask for the surfaces, for which these nets have specific properties.

**2.** Let $\bar{x}(u,v)$ be a parametric representation of a surface $\mathfrak{F} \in C^r$, $r \geq 4$, in $E^3$. We consider a moving frame $\mathfrak{D} = \{\bar{e}_i(u, v) \mid i = 1, 2, 3\}$ of $\mathfrak{F}$, such that $\bar{e}_1, \bar{e}_2$ are the principle directions, $\bar{e}_3$ is the unit normal vector and $\det(\bar{e}_1, \bar{e}_2, \bar{e}_3) = 1$.

We assume that $\mathfrak{F}$

(a) has no parabolic points and

(b) is not a canal surface.

The second assumption implies that the central surfaces are not degenerated into curves [2, p. 191].

It is known that there exist linear differential forms $\omega_1, \omega_2, \omega_{31}, \omega_{32}$ and $\omega_{12}$, so that [2]

(1) $$d\bar{x} = \omega_1 \bar{e}_1 + \omega_2 \bar{e}_2, \quad d\bar{e}_i = \omega_{i1} \bar{e}_1 + \omega_{i2} \bar{e}_2 + \omega_{i3} \bar{e}_3, \quad i = 1, 2, 3.$$

Our choice of the frame $\mathfrak{D}$ allows us to set

(2) $$\omega_{31} = -k_1 \omega_1, \quad \omega_{32} = -k_2 \omega_2, \quad \omega_{12} = q_1 \omega_1 + q_2 \omega_2,$$

where $q_1, q_2$ denote the geodesic curvatures of the lines of curvature of $\mathfrak{F}$.

We denote by $\nabla_1 f, \nabla_2 f$ the Pfaffian derivatives of a function $f(u, v) \in C^1$ along the curves

$$\omega_2 = 0, \quad \omega_1 = 0$$

respectively. Then the Mainardi-Codazzi equations can be written as follows [2, p. 175]

(3) $$\nabla_2 k_1 = q_1 (k_1 - k_2), \quad \nabla_1 k_2 = q_2 (k_1 - k_2).$$

We mention also the Theorema Egregium

(4) $$\nabla_2 q_1 - \nabla_1 q_2 - q_1^2 - q_2^2 = k_1 k_2.$$

A parametric representation of the central surface $\mathfrak{F}^{(i)}$, $i = 1, 2$, is the following

(5) $$\bar{y}^{(i)}(u,v) = \bar{x}(u,v) + \frac{1}{k_i(u,v)} \bar{e}_3(u,v).$$

We can use $\mathfrak{D}^{(1)} := \{\bar{e}_2, \bar{e}_3, \bar{e}_1\}$ resp. $\mathfrak{D}^{(2)} := \{\bar{e}_3, \bar{e}_1, \bar{e}_2\}$ as a moving frame of $\mathfrak{F}^{(1)}$ resp. $\mathfrak{F}^{(2)}$. Similarly to (1) there exist suitable linear differential forms $\omega_i^{(j)}$, $i, j = 1, 2$, such that

(6) $$d\bar{y}^{(1)} = \omega_1^{(1)}\bar{e}_2 + \omega_2^{(1)}\bar{e}_3, \qquad d\bar{y}^{(2)} = \omega_1^{(2)}\bar{e}_3 + \omega_2^{(2)}\bar{e}_1.$$

By differentiating equations (5) and making use of (1) and (6) we obtain

(7) $$\omega_1^{(1)} = \left(1 - \frac{k_2}{k_1}\right)\omega_2, \quad \omega_2^{(1)} = \frac{-dk_1}{k_1^2}, \quad \omega_1^{(2)} = \frac{-dk_2}{k_2^2}, \quad \omega_2^{(2)} = \left(1 - \frac{k_1}{k_2}\right)\omega_1.$$

We underline here that because of our assumption (b) we have

$$\omega_1^{(1)} \wedge \omega_2^{(1)} \neq 0, \quad \omega_1^{(2)} \wedge \omega_2^{(2)} \neq 0,$$

where by "$\wedge$" the wedge product of differential forms is denoted. These give, in view of (7),

(8) $$\nabla_1 k_1 \neq 0, \quad \nabla_2 k_2 \neq 0.$$

The differentials of the moving frames $\mathfrak{D}^{(1)}$, $\mathfrak{D}^{(2)}$ can be written

(9) $$d\bar{e}_2 = \omega_{12}^{(1)}\bar{e}_3 - \omega_{31}^{(1)}\bar{e}_1, \quad d\bar{e}_3 = -\omega_{12}^{(1)}\bar{e}_2 - \omega_{32}^{(1)}\bar{e}_1, \quad d\bar{e}_1 = \omega_{31}^{(1)}\bar{e}_2 + \omega_{32}^{(1)}\bar{e}_3,$$
$$d\bar{e}_3 = \omega_{12}^{(2)}\bar{e}_1 - \omega_{31}^{(2)}\bar{e}_2, \quad d\bar{e}_1 = -\omega_{12}^{(2)}\bar{e}_3 - \omega_{32}^{(2)}\bar{e}_2, \quad d\bar{e}_2 = \omega_{31}^{(2)}\bar{e}_3 + \omega_{32}^{(2)}\bar{e}_1.$$

From equations (1) and (9) we find

(10) $$\omega_{31}^{(1)} = \omega_{12}, \quad \omega_{32}^{(1)} = -\omega_{31}, \quad \omega_{12}^{(1)} = -\omega_{32},$$
$$\omega_{31}^{(2)} = -\omega_{32}, \quad \omega_{32}^{(2)} = -\omega_{12}, \quad \omega_{12}^{(2)} = \omega_{31}.$$

In addition there exist functions $a^{(i)}, b^{(i)}, c^{(i)}, q_1^{(i)}, q_2^{(i)}$, $i = 1, 2$, such that

$$\omega_{31}^{(i)} = -a^{(i)}\omega_1^{(i)} - b^{(i)}\omega_2^{(i)}, \quad \omega_{32}^{(i)} = -b^{(i)}\omega_1^{(i)} - c^{(i)}\omega_2^{(i)}, \quad \omega_{12}^{(i)} = q_1^{(i)}\omega_1^{(i)} + q_2^{(i)}\omega_2^{(i)}.$$

From these relations and (2), (7), (10) we obtain

(11) $$a^{(1)} = \frac{k_1(q_1\nabla_2 k_1 - q_2\nabla_1 k_1)}{(k_1 - k_2)\nabla_1 k_1}, \quad b^{(1)} = \frac{q_1 k_1^2}{\nabla_1 k_1}, \quad c^{(1)} = \frac{k_1^3}{\nabla_1 k_1}, \quad q_1^{(1)} = \frac{k_1 k_2}{k_1 - k_2}, \quad q_2^{(1)} = 0,$$
$$a^{(2)} = \frac{k_2^3}{\nabla_2 k_2}, \quad b^{(2)} = \frac{-q_2 k_2^2}{\nabla_2 k_2}, \quad c^{(2)} = \frac{k_2(q_2\nabla_1 k_2 - q_1\nabla_2 k_2)}{(k_1 - k_2)\nabla_2 k_2}, \quad q_1^{(2)} = 0, \quad q_2^{(2)} = \frac{k_1 k_2}{k_1 - k_2}.$$

**3.** We can now prove the following

**Proposition 1.** *The null-curves of the forms $\omega_i^{(j)}$, $i, j = 1, 2$, form an isothermic net on $\mathfrak{F}^{(j)}$ if and only if the surface $\mathfrak{F}$ is a W-surface.*

**Proof.** We show the proposition for the curves $\omega_1^{(1)} = 0$, $\omega_2^{(1)} = 0$ on $\mathfrak{F}^{(1)}$. Let $\nabla_1^{(1)} f$, $\nabla_2^{(1)} f$ be the Pfaffian derivatives of $f(u,v)$ along the curves $\omega_2^{(1)} = 0$, $\omega_1^{(1)} = 0$ respectively. From the relations

$$df = \nabla_1 f \, \omega_1 + \nabla_2 f \, \omega_2 = \nabla_1^{(1)} f \, \omega_1^{(1)} + \nabla_2^{(1)} f \, \omega_2^{(1)}$$

and (7) we have

(12) $$\nabla_1^{(1)} f = \frac{k_1 (\nabla_1 k_1 \nabla_2 f - \nabla_2 k_1 \nabla_1 f)}{(k_1 - k_2) \nabla_1 k_1}, \qquad \nabla_2^{(1)} f = \frac{-k_1^2}{\nabla_1 k_1} \nabla_1 f.$$

From (11) and (12) we obtain

$$\nabla_1^{(1)} q_1^{(1)} + \nabla_2^{(1)} q_2^{(1)} = \frac{k_1^2}{(k_1 - k_2)^3 \nabla_1 k_1} \begin{vmatrix} \nabla_1 k_1 & \nabla_2 k_1 \\ \nabla_1 k_2 & \nabla_2 k_2 \end{vmatrix}.$$

The vanishing of $\nabla_1^{(1)} q_1^{(1)} + \nabla_2^{(1)} q_2^{(1)}$ is a necessary and sufficient condition for the net

$$\omega_1^{(1)} = 0, \quad \omega_2^{(1)} = 0$$

to be isothermic while the vanishing of the determinant of the right-hand side is a necessary and sufficient condition for the surface $\mathfrak{F}$ to be a W-surface. The proposition follows then at once.

**Corollary 2.** *If one of the next properties is valid*:

 (a) *the null-curves of the forms $\omega_i^{(1)}$, $i = 1, 2$, form an isothermic net on $\mathfrak{F}^{(1)}$*,

 (b) *the null-curves of the forms $\omega_i^{(2)}$, $i = 1, 2$, form an isothermic net on $\mathfrak{F}^{(2)}$*,

 (c) *the surface $\mathfrak{F}$ is a W-surface*,

*then the rest two are also valid.*

**4.** We consider the asymptotic lines of the central surfaces $\mathfrak{F}^{(i)}$, $i = 1, 2$. These are defined by the differential equations

$$a^{(i)} \omega_1^{(i)2} + 2 b^{(i)} \omega_1^{(i)} \omega_2^{(i)} + c^{(i)} \omega_2^{(i)2} = 0, \quad i = 1, 2.$$

Taking account of (7) and (11) we easily find that to these nets correspond on the surface $\mathfrak{F}$ the nets

(13) $$\nabla_1 k_1\, \omega_1{}^2 - \nabla_1 k_2\, \omega_2{}^2 = 0,$$

(14) $$\nabla_2 k_1\, \omega_1{}^2 - \nabla_2 k_2\, \omega_2{}^2 = 0.$$

- These nets are both orthogonal if and only if

$$\nabla_1 k_1 - \nabla_1 k_2 = 0, \quad \nabla_2 k_1 - \nabla_2 k_2 = 0,$$

i.e., if and only if the difference $k_1 - k_2$ of the principle curvatures of $\tilde{\mathfrak{F}}$ is constant.

- Both nets (13) and (14) are conjugate if and only if

$$k_2 \nabla_1 k_1 - k_1 \nabla_1 k_2 = 0, \quad k_2 \nabla_2 k_1 - k_1 \nabla_2 k_2 = 0,$$

i.e., if and only if the quotient $\dfrac{k_1}{k_2}$ of the principle curvatures of $\tilde{\mathfrak{F}}$ is constant.

From the above we conclude

**Proposition 3**. *For a surface $\tilde{\mathfrak{F}}$ the nets* (13) *and* (14) *of* $\tilde{\mathfrak{F}}$ *corresponding to the asymptotic lines of the central surfaces are both*

(a) *orthogonal if and only if $\tilde{\mathfrak{F}}$ is a W-surface satisfying $k_1 - k_2 =$ const. and*

(b) *conjugate if and only if $\tilde{\mathfrak{F}}$ is a W-surface satisfying $\dfrac{k_1}{k_2} =$ const.*

**Remark.** (a) For a W-surface satisfying

$$k_1 - k_2 = \text{const.}$$

the nets (13) and (14) corresponding to the asymptotic lines of $\tilde{\mathfrak{F}}^{(i)}$, $i = 1, 2$, coincide, are real and bisect the angle between the lines of curvature of $\tilde{\mathfrak{F}}$.

(b) For a *W*-surface satisfying

$$\frac{k_1}{k_2} = \text{const.}$$

the nets (13) and (14) coincide and they are real only in case that the given surface $\tilde{\mathfrak{F}}$ is elliptic.

In view of (2) it is clear that the spherical images of the nets (13) and (14) are the following

(15) $$\frac{\nabla_1 k_1}{k_1{}^2}\, \omega_{31}{}^2 - \frac{\nabla_1 k_2}{k_2{}^2}\, \omega_{32}{}^2 = 0,$$

(16) $$\frac{\nabla_2 k_1}{k_1^2} \omega_{31}^2 - \frac{\nabla_2 k_2}{k_2^2} \omega_{32}^2 = 0.$$

These nets are both orthogonal if and only if

$$\frac{\nabla_1 k_1}{k_1^2} - \frac{\nabla_1 k_2}{k_2^2} = 0, \quad \frac{\nabla_2 k_1}{k_1^2} - \frac{\nabla_2 k_2}{k_2^2} = 0,$$

i.e., if and only if the difference

$$\frac{1}{k_1} - \frac{1}{k_2}$$

of the principle radii of curvature is constant. Thus we have

**Proposition 4.** *The spherical images* (15) *and* (16) *of the nets* (13) *and* (14) *corresponding to the asymptotic lines of the central surfaces of a surface* $\mathfrak{F}$ *are both orthogonal if and only if* $\mathfrak{F}$ *is a W-surface satisfying*

$$\frac{1}{k_1} - \frac{1}{k_2} = const.$$

**5.** We concentrate now our point of interest on the lines of curvature of the central surfaces $\mathfrak{F}^{(i)}$, $i = 1, 2$. These are defined by the differential equations

$$b^{(i)} \omega_1^{(i)2} + \left(c^{(i)} - a^{(i)}\right) \omega_1^{(i)} \omega_2^{(i)} - b^{(i)} \omega_2^{(i)2} = 0, \quad i = 1, 2.$$

By virtue of (7) and (11) to these nets correspond on the surface $\mathfrak{F}$ the nets

(17) $$q_1 \nabla_1 k_1 \omega_1^2 + \left[k_1^2(k_1 - k_2) + q_2 \nabla_1 k_1 + q_1 \nabla_2 k_1\right] \omega_1 \omega_2 + q_2 \nabla_2 k_1 \omega_2^2 = 0,$$

(18) $$q_1 \nabla_1 k_2 \omega_1^2 + \left[k_2^2(k_1 - k_2) + q_2 \nabla_1 k_2 + q_1 \nabla_2 k_2\right] \omega_1 \omega_2 + q_2 \nabla_2 k_2 \omega_2^2 = 0.$$

These nets are both orthogonal if and only if

(19) $$q_1 \nabla_1 k_1 + q_2 \nabla_2 k_1 = 0, \quad q_1 \nabla_1 k_2 + q_2 \nabla_2 k_2 = 0.$$

Using the Mainardi-Codazzi equations (3) the preceding relations become

(20) $$q_1 \nabla_1 k_1 = -q_1 q_2 (k_1 - k_2), \quad q_2 \nabla_2 k_2 = -q_1 q_2 (k_1 - k_2).$$

We distinguish two cases:

I. The surface $\mathfrak{F}$ is a *moulding surface, i.e.,* $\mathfrak{F}$ is generated by the orthogonal trajectories of a one-parameter family of planes. Then $q_1$ or $q_2$ vanishes [1, p. 192]. Let $q_1$ vanish. From relations (20) we obtain

$$q_2 \nabla_2 k_2 = 0.$$

But $q_2 \neq 0$ because of (4) and our assumption (a). Consequently $\nabla_2 k_2 = 0$, which contradicts (8).

II. The surface $\mathfrak{F}$ is not a moulding surface. Then (20) reduce to

$$\nabla_1 k_1 = -q_2 (k_1 - k_2), \quad \nabla_2 k_2 = -q_1 (k_1 - k_2).$$

Using again (3) we obtain $\nabla_1(k_1 + k_2) = \nabla_2(k_1 + k_2) = 0$, i.e., $d(k_1 + k_2) = 0$, which means that $\mathfrak{F}$ is a surface of constant mean curvature.

Conversely, if $\mathfrak{F}$ is a surface of constant mean curvature, then it easily follows that (19) hold true and the nets (17), (18) are orthogonal.

Let now suppose that both nets (17) and (18) are conjugate. Then we have

(21) $\qquad k_2 q_1 \nabla_1 k_1 + k_1 q_2 \nabla_2 k_1 = 0, \quad k_2 q_1 \nabla_1 k_2 + k_1 q_2 \nabla_2 k_2 = 0.$

Relations (21) are impossible in the case that $\mathfrak{F}$ is a moulding surface. If $\mathfrak{F}$ is not a moulding surface relations (21) become in virtue of (3)

$$k_2 \nabla_1 k_1 = -k_1 q_2 (k_1 - k_2), \quad k_1 \nabla_2 k_2 = -k_2 q_1 (k_1 - k_2).$$

Then we obtain $\nabla_1(k_1 k_2) = \nabla_2(k_1 k_2) = 0$, i.e., $d(k_1 k_2) = 0$, which means that $\mathfrak{F}$ is a surface of constant Gaussian curvature.

Conversely, if $\mathfrak{F}$ is a surface of constant Gaussian curvature, then we can verify that (21) are valid. Hence we obtain the following

**Proposition 5**. *For a surface $\mathfrak{F}$, which is not a moulding surface, the nets* (17) *and* (18) *of $\mathfrak{F}$ corresponding to the lines of curvature of the central surfaces are both*

(a) *orthogonal if and only if $\mathfrak{F}$ is a surface of constant mean curvature and*

(b) *conjugate if and only if $\mathfrak{F}$ is a surface of constant Gaussian curvature.*

Finally we prove the following

**Proposition 6**. *Under the hypothesis of Proposition 5, the spherical images of the nets* (17) *and* (18) *are both orthogonal if and only if $\mathfrak{F}$ is a W-surface satisfying*

$$\frac{1}{k_1}+\frac{1}{k_2} = const.$$

**Proof.** Using (2) we find that the spherical images of the nets (17) and (18) are defined by

$$q_1 k_2^2 \nabla_1 k_1 \omega_{31}^2 + k_1 k_2 \left[ k_1^2 (k_1 - k_2) + q_2 \nabla_1 k_1 + q_1 \nabla_2 k_1 \right] \omega_{31}\omega_{32} + q_2 k_1^2 \nabla_2 k_1 \omega_{32}^2 = 0,$$

$$q_1 k_2^2 \nabla_1 k_2 \omega_{31}^2 + k_1 k_2 \left[ k_2^2 (k_1 - k_2) + q_2 \nabla_1 k_2 + q_1 \nabla_2 k_2 \right] \omega_{31}\omega_{32} + q_2 k_1^2 \nabla_2 k_2 \omega_{32}^2 = 0.$$

These nets are both orthogonal if and only if

$$q_1 k_2^2 \nabla_1 k_1 + q_2 k_1^2 \nabla_2 k_1 = 0, \quad q_1 k_2^2 \nabla_1 k_2 + q_2 k_1^2 \nabla_2 k_2 = 0.$$

We exclude again the moulding surfaces. The preceding relations combined with the Mainardi-Codazzi equation (3) yield

$$\nabla_1 k_2 = -\frac{k_2^2}{k_1^2} \nabla_1 k_1, \quad \nabla_2 k_2 = -\frac{k_2^2}{k_1^2} \nabla_2 k_1,$$

therefore we have $d(\frac{1}{k_1}+\frac{1}{k_2}) = 0$, i.e. the sum of the principle radii of curvature is constant.

Stylianos Stamatakis
Department of Mathematics
Aristotle University of Thessaloniki
54006 Thessaloniki, Greece
e-mail: stamata@math.auth.gr